\documentclass[12pt,twoside]{amsart}
\usepackage{amssymb}
\usepackage{amscd}

\newcommand{\kk}{\mathbf{k}}
\newcommand{\zk}{\zeta({\kk})}
\newcommand{\un}{\mathcal{U}n}
\newcommand{\uns}{\un(X)}
\newcommand{\unt}{\un(T)}
\newcommand{\Spec}{\operatorname{Spec}}

\newcommand{\Z}{\mathbf{Z}}
\newcommand{\vv}{\mathcal{V}}
\newcommand{\spec}{\operatorname{Spec}}
\newcommand{\Li}{\operatorname{Li}}
\newcommand{\End}{\operatorname{End}}
\newcommand{\gr}{\operatorname{Gr}}
\newcommand{\res}{\operatorname{Res}}
\newcommand{\tbar}{\bar{t}}
\newcommand{\dlog}{\operatorname{dlog}}
\newcommand{\hr}{H_{\textup{rig}}}
\newcommand{\unad}{Un(\Adag,D)}
\newcommand{\unx}{Un(X,D)}

\newcommand{\tcat}{\mathcal{T}}
\newcommand{\scat}{\mathcal{S}}
\newcommand{\vecs}{\operatorname{Vec}}
\newcommand{\rep}{\operatorname{Rep}}
\newcommand{\aut}{\operatorname{Aut}}
\newcommand{\autt}{\aut^\otimes}
\newcommand{\Homm}{\operatorname{Hom}^\otimes}
\newcommand{\Q}{\bold Q}
\newcommand{\cX}{\mathcal{X}}
\newcommand{\cY}{\mathcal{Y}}
\newcommand{\Hom}{\operatorname{Hom}}
\renewcommand{\O}{\mathcal{O}}
\newcommand{\zz}{{00}}
\newcommand{\Fr}{\operatorname{Fr}}

\newcommand{\Adag}{A^\dagger}
\newcommand{\Bdag}{B^\dagger}
\newcommand{\dr}{\mathcal{DR}}
\newcommand{\locs}{\operatorname{Loc}}

\title[The Double Shuffle Relations]
{The Double Shuffle Relations for $p$-adic Multiple Zeta Values}
\author{Amnon Besser}
\address{
Department of Mathematics\\
Ben-Gurion University of the Negev\\
P.O.B. 653\\
Be'er-Sheva 84105\\
Israel
}
\email{bessera@math.bgu.ac.il}
\thanks{The first author was supported in part by an Israel Science Foundation grant}
\author{Hidekazu Furusho}
\address{
Graduate School of Mathematics\\ 
Nagoya University\\
Chikusa-ku\\
Furo-cho\\
Nagoya\\
464-8602\\
Japan 
}
\email{furusho@math.nagoya-u.ac.jp}

\newtheorem{thm}{Theorem}[section]
\newtheorem{lem}[thm]{Lemma}
\newtheorem{cor}[thm]{Corollary}
\newtheorem{prop}[thm]{Proposition}  

\theoremstyle{remark}

\theoremstyle{definition}
\newtheorem{defn}[thm]{Definition}
\newtheorem{rem}[thm]{Remark}
\newtheorem{nota}[thm]{Notation}

\newtheorem{pf}{Proof}

\numberwithin{equation}{section}

\begin{document}
\maketitle

\begin{abstract}
We give a proof of double shuffle relations for $p$-adic multiple zeta values
by developing higher dimensional version of tangential base points and 
discussing a relationship with two (and one) variable $p$-adic multiple
polylogarithms. 
\end{abstract}

\tableofcontents
\setcounter{section}{-1}
\section{Introduction}
\label{sec:intro}

In this paper we will prove a set of formulas, known as double
shuffle relations, relating the $p$-adic multiple zeta values defined
by the second named author in~\cite{F}. These formulas are analogues
of formulas for the usual (complex) multiple zeta values. These have a
very simple proof which unfortunately does not translate directly to
the $p$-adic world.

Recall that the (complex) multiple zeta value $\zk$, where $\kk$
stands for the multi-index $\kk=(k_1,\ldots, k_m) $, is defined by the formula
\begin{equation}\label{defform}
 \zk = \underset{n_i\in\bold  N}{\underset{0<n_1<\dotsm<n_m}{\sum}}
 \frac{1}{n_1^{k_1}\cdots n_m^{k_m}}\;,
\end{equation}
The series is easily seen to be convergent assuming that $k_m > 1$.

Multiple zeta values satisfy two types of so called shuffle product formulas,
expressing a product of multiple zeta values as a linear combination
of other such values. The first type of formulas are called series
shuffle product  
formulas (sometimes called by harmonic product formulas). 
The simplest example is the relation 
\begin{equation}\label{serrel}
\zeta(k_1)\cdot \zeta(k_2)= \zeta(k_1,k_2) + \zeta(k_2,k_1) +
\zeta(k_1+k_2)\;,
\end{equation}
which is easily obtained from the expression \eqref{defform}
by noting that the left hand side is a summation over an infinite
square of pairs $(n_1,n_2)$
of the summand in \eqref{defform},
and that summing over the lower triangle (respectively the upper
triangle, respectively the diagonal) gives the three terms on the
right hand side. Every series shuffle product formula has this type of
proof.

The second type of shuffle product formulas, known as iterated
integral shuffle product formulas, is somewhat harder to
establish and follows from the description of multiple zeta values in
terms of multiple polylogarithms. More precisely. The one variable multiple
polylogarithm is defined by the formula

\begin{equation}\label{multpolll}
  \Li_\kk(z) = \underset{n_i\in\bold  N}{\underset{0<n_1<\dotsm<n_m}{\sum}} 
  \frac{z^{n_m}}{n_1^{k_1}\cdots n_m^{k_m}}\;,
\end{equation}
near $z=0$. It can then be extended as a multi-valued function to
$\bold P^1(\bold C)-\{0,1,\infty\}$. We clearly have the relation 
$\lim_{z\to 1} \Li_\kk(z)=\zk$.

Multiple polylogarithms can be written using the theory of iterated
integrals due to Chen~\cite{Ch}. In other words, they satisfy
a system of unipotent differential equations. This gives an integral
expression for multiple polylogarithms. By substituting $z=1$ and
splitting the domain of integration in the right way we obtain the
iterated integral shuffle product formulas, a simple example of which is
the formula
\begin{equation}\label{intshuff}
\zeta(k_1)\cdot\zeta(k_2)=\sum_{i=0}^{k_1-1}\binom{k_2-1+i}{i}\zeta(k_1-i,k_2+i)+\sum_{j=0}^{k_2-1}\binom{k_1-1+j}{j}\zeta(k_2-j,k_1+j).
\end{equation}

In~\cite{F} the second named author defined the $p$-adic version of
multiple zeta values and studied some of their properties. The
defining formula \eqref{defform} can not be directly used $p$-adically
because the defining series does not converge. Instead, one must use
an indirect approach based on the theory of Coleman
integration~\cite{Co,Be}. Coleman's theory defines $p$-adic analytic
continuation for solutions of unipotent differential equations ``along
Frobenius''. Coleman used his theory initially to define $p$-adic
polylogarithms. In~\cite{F} Coleman integration was used to define one
variable $p$-adic multiple polylogarithms. Taking the limit at $1$ in
the right way one
obtains the definition of $p$-adic multiple zeta values. It is by no
means trivial that the limit even exists or is independent of choices,
and this is the main result of~\cite{F}.

Given their definition, it is not surprising that for $p$-adic
multiple zeta values it is the iterated
integral shuffle product formulas that are easier to obtain
$p$-adically. In~\cite{F} the series shuffle product formulas
were not obtained. The purpose of this work is to prove
(Theorem~\ref{main}) these formulas, and as a consequence the double
shuffle relations (Corollary~\ref{maincor}) for $p$-adic multiple zeta values.

To prove the main theorem it is necessary to use the theory of
Coleman integration in several variables developed by the first named author
in~\cite{Be}. The reason for this is quite simple - If one tries to
replace multiple zeta values by multiple polylogarithms in the proof
of \eqref{serrel} sketched above one easily establishes the formula
\begin{equation}\label{tvarfor}
  \Li_{k_1}(z) \Li_{k_2}(w)= \Li_{k_1,k_2}(z,w) + \Li_{k_2,k_1}(w,z)
  + \Li_{k_1+k_2}(zw)\;,
\end{equation}
which is a two variable formula. It seems impossible to obtain a one
variable version of the same formula. The proof of the main theorem
thus consists roughly speaking of showing that \eqref{tvarfor} extends
to Coleman
functions of several variables and then taking the limit at $(1,1)$.

Since taking the limit turned out to be rather involved in~\cite{F},
we opted for an alternative approach, which was motivated by a letter
of Deligne to the second named author~\cite{De2}. Deligne observes
that taking the limit at $1$ for the multiple polylogarithm can be
interpreted as doing analytic continuation from tangent vectors at $0$
and $1$, using the theory of the tangential basepoint at infinity
introduced in~\cite{De1}. To analytically continue \eqref{tvarfor}
and obtain the series shuffle product formula we analyze a more general
notion of tangential basepoint sketched in loc.\ cit.\ and examine
among other things its
relation with Coleman integration.

To give a precise meaning of the limit value to $(1,1)$,
we work over the moduli space $\mathcal{M}_{0,5}$ of curves of 
$(0,5)$-type and the normal bundles for the divisors at infinity
$\overline{\mathcal{M}_{0,5}}-\mathcal{M}_{0,5}$
($\overline{\mathcal{M}_{0,5}}$: a compactification of $\mathcal{M}_{0,5}$).
Two variable $p$-adic multiple polylogarithms are introduced.
They are Coleman functions over $\mathcal{M}_{0,5}$.
In \S \ref{anacon}
their analytic continuation to the normal bundle  will be discussed.
In particular, we will relate the behavior of the analytic continuation of 
two variable multiple polylogarithm to a normal bundle 
with one variable multiple polylogarithm and then
we get $p$-adic multiple zeta values as ``special values''
of two variable multiple polylogarithms.

It is a pleasure to thank P. Deligne and A. Goncharov. 
The first named author was
supported by a grant from the Israel Science foundation. The second
named author is supported in part by
JSPS Research Fellowships for Young Scientists. He would like to thank
the Center for Advanced Studies in
Mathematics at Ben-Gurion University for supporting a visit to
Beer-Sheva in February 2003 when a large part of the research
described in this work was carried out and he would also like to thank
the Johns Hopkins University for their support and hospitality during
the JAMI conference.


\section{Tannakian duality and fundamental groups}
\label{sec:tannak}

This section is introductory, and contains some basics of the theory
of Tannakian categories (for which the standard reference
is~\cite{DM}) and some ways it can be used to construct fundamental
groups and path spaces in algebraic geometry.

The theory of Tannakian categories has its roots in Tannaka-Krein
duality, itself a generalization of Pontrjagin duality, which shows
how to reconstruct a compact topological group out of the category of
its representations. Subsequently, it was realized by Grothendieck
that the reconstruction can be done, in the somewhat different context
of affine group schemes, based on some formal properties of this
category of representations. Consequently, any category with some
additional structures and properties, formalized in the notion of a
Tannakian category, gives rise to an affine group scheme.

A Tannakian category (or more precisely a neutral Tannakian category)
over a field $F$ is a rigid abelian $F$-linear tensor category
possessing a fiber functor
\cite[Definition~2.19]{DM}. An abelian $F$-linear category $\tcat$ becomes a
tensor category with the addition of a bifunctor $\otimes :\tcat\times
\tcat \to \tcat$ satisfying certain axioms mimicking those of the
tensor product of vector spaces. It is rigid if it satisfies
\cite[Definition 1.7]{DM}. A tensor functor $G: \tcat \to \tcat'$ is a
functor together with a functorial isomorphism $G(X)\otimes G(Y) \to
G(X\otimes Y) $ satisfying a few obvious properties \cite[Definition~1.8]{DM}

The category $\vecs_F$ of
finite dimensional vector spaces over $F$ is a rigid abelian
$F$-linear tensor category. The final data in the
definition of a Tannakian category is a fiber functor, which by
definition is an exact $F$-linear faithful tensor functor $\omega: \tcat
\to \vecs_F $. 

Suppose now that $G$ is an affine group scheme over $F$. One than has
the category $\rep_F(G)$ of finite dimensional $F$-rational
representations of $G$. This is clearly a Tannakian category with
fiber functor sending a representation to its underlying $K$-vector
space. The main result of Tannakian duality \cite[Theorem~2.11]{DM}
says that $\rep_F(G)$ 
determines $G$ uniquely, and conversely, any Tannakian category
arises as $\rep_F(G)$ for some $G$.

Recall that an affine group scheme over $F$ can be defined as a group
valued functor on the category of $F$-algebras which is represented by
an $F$-algebra. The way that $\rep_F(G)$ determines $G$ is as follows:
$G$ is isomorphic to the group valued functor $\autt(\omega) $ whose
value on an $F$-algebra $R$ is
\begin{equation*}
  \autt(\omega)(R) := \Hom(\omega\otimes R,\omega\otimes R)\;,
\end{equation*}
where $\omega\otimes R $ is the composition of $\omega$ with the
change of ring functor and $\Hom$ is in the category of tensor
functors. We will often write $\pi_1(\tcat,\omega) $ instead of $
\autt(\omega) $ to indicate the dependency on the category $\tcat$.

Tannakian duality can be used to define all kinds of ``fundamental
groups'' in algebraic geometry and algebraic topology.
The basic example to keep in mind is that of the category $\locs_X$ of
local system 
of finite dimensional $\Q $-vector spaces on a connected topological
space $X$. Given a 
choice of a point $x\in X$, this is equivalent to the category of
finite dimensional $\Q$-representations of the fundamental group
$\pi_1(X,x) $. It is clearly a Tannakian category over $\Q$. Each
point $x\in X$ determines a fiber 
functor $\omega_x$ where $\omega_x(V)=$ fiber of $V$ at $x$. The
obvious map $\pi_1(X,x) \to \pi_1(\locs_X,\omega_x) (\Q) $ is
explicitly given by
\begin{equation}\label{mon1}
  \gamma \mapsto \left(V\mapsto \text{monodromy along the loop }
    \gamma \in \aut(\omega_x(V))\right) \;.
\end{equation}

We may also use Tannakian duality to define generalized path
spaces. Namely, if $\tcat$ is a Tannakian category and $\omega_1$ and
$\omega_2$ are two fiber functors on $\tcat$. Then there is a affine
scheme $\pi_1(\tcat,\omega_1,\omega_2)=\Homm(\omega_1,\omega_2)$ whose
points in an $F$-algebra $R$ are
\begin{equation*}
  \Homm(\omega_1,\omega_2)(R) := \Hom(\omega_1\otimes
  R,\omega_2\otimes R)\;,
\end{equation*}
Again Hom taken in the category of tensor functors.
There is an obvious composition map,
\begin{equation}
  \label{eq:picompos}
  \circ : \pi_1(\tcat,\omega_2,\omega_3) \times
  \pi_1(\tcat,\omega_1,\omega_2) \to \pi_1(\tcat,\omega_1,\omega_3)\;,
\end{equation}
and in particular $\pi_1(\tcat,\omega_1,\omega_2)$ is a homogeneous
space for an action from the 
left of $\pi_1(\tcat,\omega_1) $ as well as for an action from the
right of $\pi_1(\tcat,\omega_2) $ and is principal for both
actions. For $\tcat = \locs_X$ we have a map $\pi_1(X,x,y) \to
\pi_1(\locs_X,\omega_x,\omega_y)(\Q) $, where $\pi_1(X,x,y)  $ is the
space of paths between $x$ and $y$ on $X$ up to homotopy, given by
\begin{equation}
  \label{mon2}
    \gamma \mapsto \left(V\mapsto \text{parallel transport along the path }
    \gamma \in \Hom(\omega_x(V),\omega_y(V))\right) \;.
\end{equation}

A tensor functor $F:\scat \to \tcat $ clearly induces a map of path
spaces
\begin{equation}
  \label{eq:pifunc}
  F^\ast :  \pi_1(\tcat,\omega_1,\omega_2) \to
  \pi_1(\scat,\omega_1\circ F,\omega_2\circ F)
\end{equation}
which is explicitly given by sending $\gamma \in
\Homm(\omega_1\otimes R,\omega_2\otimes R) $ to $F^\ast (\gamma)$ given
by $ F^\ast(\gamma)_S = \gamma_{F(S)} $. This is clearly compatible
with composition. In the case of categories of local systems, when
$F$ is induced by the 
continuous map $f: X \to Y$ it is easy to see that this map is
compatible, via \eqref{mon2} with the usual map induced by $f$ on path spaces.

\section{Coleman's $p$-adic integration}
\label{sec:coleman}

In this section we recall the theory of Coleman's $p$-adic integration in
several variables introduced in~\cite{Be} and make some additional
constructions which are required for the current work. In some
respects, the similar theory of Vologodsky~\cite{V} is more suitable
for the current setting. However, because it is much less explicit it
seems harder to adapt it to the tangential constructions of
\S \ref{sec:tangential}.

Recall from~\cite{Be} that the basic setup for Coleman integration
theory is as follows: We have a field of characteristic $p$, $\kappa$,
which we assume for simplicity to be algebraically closed and a discrete
valuation ring
$\vv$ with fraction field $K$ and residue field $\kappa$.
A rigid triple is a triple $T=(X,Y,P)$ consisting of $P$, a formal $p$-adic
$\vv$-scheme, $Y$ a closed $\kappa$-subscheme of $P$ which is proper
over $\spec(\kappa)$ and $X$ an open $\kappa$-subscheme of $Y$ such
that $P$ is smooth in a neighborhood of $X$.

A simple case where rigid triples arise
(see~\cite[Definition~5.1]{Be}) is that of tight rigid triples.
Let $\cX\subset \cY$ be an open immersion of $\vv$-schemes such that
$\cX$ is smooth and $\cY$ is complete. The associated triple
$T_{(\cX,\cY)}:=(\cX\otimes_\vv
\kappa,\cY\otimes_\vv \kappa,\hat{\cY})$, where $\hat{\cY}$ is the
$p$-adic completion of $\cY$ is called a tight rigid triple. An affine
rigid triple is a tight rigid triple $T_{(\cX,\cY)}$ with $\cX$ affine.

Given a rigid triple as above we have the category of unipotent overconvergent
isocrystals on $T$, denoted $\unt$. This is the category of sheaves on
the tube of $X$ in $P$, overconverging into the tube of $Y$ in $P$ in
the sense of Berthelot \cite{Ber}, with an integrable connection. This category
depends, up to a unique isomorphism, only on $X$, so we will often
denote it simply by $\uns$.

A Frobenius endomorphism of $X$ is a $\kappa$-linear endomorphism of
$X$, which is some power of the geometric Frobenius. A Frobenius
endomorphism $\phi$ 
induces an auto-functor $\phi^\ast$ on $\uns$. Suppose that
$T=T_{(\cX,\cY)}$ is a tight rigid triple and that $\phi $ is the
reduction of an endomorphism $\varphi$ of the pair $(\cX,\cY)$. Then
the functor  $\phi^\ast$ is given by the obvious functor induced by
$\varphi$ on $\unt$.
\begin{defn}\label{frobfunc}
A Frobenius fiber functor on $\uns$ is a fiber functor $\omega$
together with an isomorphism of fiber functors (which we write as an
equality to simplify matters) $\omega\circ \phi^\ast = \omega $
\end{defn}

For two Frobenius fiber functors $\omega_1$ and $\omega_2$ there is an
obvious action of $\phi^\ast$ on the path space $ \pi_1(\uns,\omega_1,\omega_2)$ coming from the action \eqref{eq:pifunc}.
The main result of~\cite{Be} (Corollary 3.2) can be stated, in a
slightly generalized form, as the following:
\begin{thm}\label{mainbe}
Given any two Frobenius
fiber functors on $\uns$, $\omega_1$ and $\omega_2$, there exists a canonical
invariant path, i.e., an isomorphism of fiber functors,
$a_{\omega_1,\omega_2}:\omega_1\to\omega_2$,
which is fixed by 
$\phi^\ast$.
We will call the isomorphism $\omega_1 \to \omega_2$ analytic
continuation along Frobenius. It has the basic compatibility property that
$a_{\omega_2,\omega_3} \circ a_{\omega_1,\omega_2} =
a_{\omega_1,\omega_3}$. 
\end{thm}
\begin{proof}
We only sketch the main idea, which is rather simple, and refer to
reader to loc.\ cit.\ for the full story. The key point is that for
any Frobenius fiber functor $\omega$ there is an action $\phi^\ast$
on the fundamental group $\pi_1(\uns)$, which is a pro-unipotent
group, and its Lie algebra is controlled by negative tensor powers of
$\hr^1(X/K) $. The action of $\phi^\ast$ on $\hr^1(X/K) $ has strictly
positive weights, and consequently one can show that the map
\begin{equation*}
  g\mapsto \phi^\ast (g) g^{-1}\;,\;  \pi_1(\uns,\omega) \to
  \pi_1(\uns,\omega)\;, 
\end{equation*}
is a bijection. Now, the space $\pi_1(\uns,\omega_1,\omega_2) $ is a
$\pi_1(\uns,\omega_1)$-principal 
space and its $\phi^\ast$ action is compatible with the group
action. From this it is very easy to deduce that there exists a unique 
$\phi^\ast$-invariant element. The uniqueness immediately implies the
compatibility.
\end{proof}

Usually, Frobenius fiber functors as above are obtained from geometric points
$x:\Spec \kappa\to X$ by pullback, $\omega_x = x^\ast$, assuming that $x$
is fixed by $\phi$. Concretely, for an overconvergent isocrystal
$(M,\nabla)$ the fiber functor $\omega_x(M,\nabla) $ is
realized as the vector space of horizontal sections of $\nabla$ on the
residue disc
\begin{equation*}
  U_x = ]x[_P\; \text{(=points reducing to $x$)}
\end{equation*}
which is the tube in the sense of Berthelot.
Since each point of $X$ will be fixed by some
Frobenius automorphism it is not hard to see that the invariant path
is independent of the choice of the Frobenius endomorphism. 

An abstract Coleman function on $T$ is defined to be a triple $(M,s,h)$ where
  \begin{itemize}
  \item $M=(M,\nabla)$ is a unipotent isocrystal on $T$.
  \item $s\in \Hom(M,\mathcal{O})$.
  \item $h$ is a collection of sections, $\{h_x\in M(U_x),\; x\in
    X\}$, with
    $\nabla(h_x)=0$, which correspond to each other via analytic
    continuation along Frobenius.
  \end{itemize}

One can evaluate an abstract Coleman function on each residue disc
$U_x$ by taking $s(h_x)$. Coleman functions are then equivalence classes of
abstract Coleman functions under a relation which essentially
guarantees that their evaluations will be identical. In fact, a
Coleman function is a connected component of the category of abstract
Coleman functions, where a morphism $(M_1,s_1,h_1) \to (M_2,s_2,h_2) $
is a horizontal map $u: M_1 \to M_2$ such that $s_2\circ u = s_1$ and
$u(h_1)=h_2$. It is
shown in~\cite[Corollary 4.13]{Be} that two Coleman functions which
coincide on
an open subset of the tubular neighborhood $]X[_P$ of $X$ in $P$ 
coincide everywhere.
\begin{rem}\label{integclos}
It will be important to note how the theory can be used to define
iterated integrals. For a similar and more detailed discussion
see the proof of Theorem~4.15 in~\cite{Be}. Suppose that we have a
collection of Coleman
functions $F_i$, arising from the abstract Coleman functions
$(M_i,s_i,h_i)$, and one-forms $\omega_i$ for $i=1,\ldots, k$. Suppose now
that the locally analytic differential form $\omega=\sum F_i \omega_i $ is
closed, and that we wish to find a Coleman function $F$ such that $dF =
\omega$. We may begin by constructing the connection $M = (\oplus M_i)
\oplus \O_T$ with the connection given by the formula
\begin{equation*}
  \nabla(m_1,\cdots,m_k,f) = (\nabla_1 m_1,\cdots ,\nabla_k m_k,
  df-\sum \omega_i s_i(m_i))\;.
\end{equation*}
This is clearly a unipotent connection. If it is integrable we
immediately obtain our required Coleman function represented by the
abstract Coleman function $(M,s,y) $, where $s$ is the projection on
the last factor and $h$ is obtained by choosing one residue disc $U$,
extending the sum of the $h_i$ to a horizontal section of $M$ on $U$,
and then analytically continuing along Frobenius. In general, $\nabla$ will
not be integrable. However, there is a maximal integrable
subconnection $(M^\textup{int},\nabla)$ of $(M,\nabla)$
\cite[(2.3)]{Be}. The fact 
that $\sum F_i \omega_i$ is closed guarantees that we may choose $h$ to
be in $M^\textup{int}(U)$ hence we obtain our Coleman function from
$(M^\textup{int},s,h) $.
\end{rem}
One can apply the theory to other fiber functors. A simple case was
already discussed in~\cite[Section~5]{Be}. Here we will need an easy
generalization (to the two dimensional case for simplicity) as
follows: Suppose that $T=(X,Y,P)=T_{(\cX,\cY)}$ is a rigid
triple and that $y\in D=Y-X$ is a closed point and that $D$ is locally
given near $y$ by the reduction of two parameters $t_1$ and $t_2$ on
$Y$. Let $A$ be the ring of Laurent series with coefficients in $K$,
\begin{equation*}
 A:=  \Bigl\{f(t_1,t_2)=\sum_{(i,j)\in \Z^2} a_{ij}t_1^it_2^j 
\;\Bigm|\; f \text{
 converges on } r < |t_1|, |t_2|<1 \text{ for some } r\Bigr\}
\end{equation*}
It follows from~\cite{Ber} that overconvergent isocrystals on $T$
give rise to a connection on a (free) $A$-module.
\begin{prop}
  An overconvergent unipotent connection on $T$ has a full set of
  solutions in $B:=A[\log(t_1),\log(t_2)]$.
\end{prop}
\begin{proof}
  Overconvergent unipotent crystals are solved by iterated integration
  and the ring $B$ is easily seen to be closed under partial
  integration with respect to $t_1$ and $t_2$.
\end{proof}
Note that the ring $B$ is independent of choice of parameters $t_1$
and $t_2$
. The following operation will be required
later.
\begin{defn}\label{contstterm}
  The constant term with respect to the parameters $t_1$ and $t_2$ of
  $f\in B$ is defined as follows: Let $f=\sum f_{ij} \log^i(t_1)
  \log^j(t_2) $. Let $f_{00}= \sum  a_{ij}t_1^it_2^j$. Then the
  constant term is $a_{00}$
\end{defn}
Note that the constant term does depend on the choice of parameters.
\begin{defn}\label{nfiber}
  The fiber functor $\omega_y$ of $\uns$ associates to a unipotent
  overconvergent isocrystal the vector space of solutions in $B$. 
\end{defn}

Suppose that we have an endomorphism $\varphi $ of the pair
$(\cX,\cY)$ reducing to a Frobenius endomorphism fixing $y$ (a
sufficiently high power of a Frobenius endomorphism will fix $y$).
Then, pulling back by $\varphi$ gives an isomorphism $\omega_y \circ
\varphi \cong \omega_y $ making it a Frobenius fiber functor in the
sens of Definition~\ref{frobfunc}. It is therefore possible to analytically
continue along Frobenius to or from $\omega_y$.

To end this section we note that by the naturality of the analytic
continuation along Frobenius, it clearly extends to pro-unipotent
isocrystals. We also note that as long as the underlying connections
are defined over a discretely valued subfield it is possible to work
over $\bold{C}_p$.

\section{Tangential basepoints}
\label{sec:tangential}

In this section we first recall the theory of the tangential basepoint
in the de Rham setting due to Deligne~\cite{De1}. In its simplest
form this theory allows
one to define fiber functors for the category of integrable
connections on a curve associated with a tangent point ``at
infinity''. We explain Deligne's interpretation of constant terms in
terms of these tangential basepoints. In loc.\ cit.\ a
higher dimensional theory is also sketched with no details, but these
can easily be filled in, as we will do. We then discuss the possibility
of analytically continuing solutions of a unipotent integrable
connection along Frobenius a la Coleman to such fiber functors. We
finally explain how one can iterate the construction of tangential
basepoints in two different ways and get the same result and we close
with some applications to Coleman integration theory.

We first recall the construction of Deligne~\cite[15. Theorie
alg\'{e}brique]{De1}. Let $C$ be a curve over a field $K$ of
characteristic $0$, smooth at a point $P$, and
let $t$ be a parameter at $P$ (Deligne immediately passes to the
completion at $P$ but we will not do the same). Suppose $(M,\nabla)$ is a
connection on $C-\{P\}$ with logarithmic singularities along $P$. This means
that locally near $P$ the connection $\nabla$ can be written as
$\nabla = d+ \Gamma$,
where $\Gamma$ is a section of $\End(M)\otimes\Omega^1_C(\log P)$ with at most a simple pole at
$P$. The parameter $t$ induces naturally a parameter $\tbar$ on the
tangent space $T_P(C)$ (the linear parameter taking the value $1$ at
the derivation $d/dt$). We associate with this data the connection on
$T_P(C)-0 $ on the trivial bundle with fiber $M_P$ (fiber of $M$ at
$P$), given by $\res_P(\nabla) := d +(\res_P \Gamma)
\dlog(\tbar)$. This clearly defines a functor $\res_P $.

An easy computation shows that the functor $\res_P$ does not depend
on the parameter $t$.  Deligne gives an
alternative, coordinate free description of the same construction,
making this fact evident. The valuation $v_P$ on the fraction field
$K(C)$ gives an algebra filtration $F_P$ on this field and there is an obvious
canonical isomorphism between the associated graded algebra $\gr_PK(C)
$ and the coordinate ring of $T_P(C)-0$, given by sending a cotangent
vector $d_P f $, thought of as a linear function on $T_P(C) $, to the
image of $f$ in $\gr_P^1 K(C) $.

Let $ \Omega_{\O_C/K}^1(\log P)$ be the sheaf of differential forms
with log singularities along $P$.
We give $\Omega_{K(C)/K}^1 = K(C) \otimes_{\O_C}
\Omega_{\O_C/K}^1(\log P)$ the filtration induced from the filtration
on the first term.
The filtrations on $K(C)$ and on $\Omega_{K(C)/K}^1 $ induces a
filtration on $M\otimes_{\O_C}K(C) $ and on $M\otimes_{\O_C}
\Omega_{K(C)/K}^1$ and the assumption that $\nabla$ has log
singularities implies that it preserves the filtration. It is easy to
see that the induced connection on the associated graded is exactly
$\res_P(\nabla) $

The construction of $\res_P $ gives us the option of producing more
fiber functors for the category of connections on $C$ by taking the
fiber at a point of $T_P(C)-0$. A particular case of this construction
gives the notion of a constant term for a horizontal section of this
connection as we now explain. Suppose that the connection $\nabla $
is unipotent. Then it is very easy to see that one can find a basis of
formal
horizontal sections to $\nabla$ near $P$ with coefficients in the ring
$K[[t]][\log(t)] $, where 
$\log(t)$ is treated as a formal variable whose derivative is
$dt/t$.  Let $v$ be such a horizontal section.  There is a sense in which we can
specialize $v$ to the fiber $M_P$, namely, taking
the constant term.
\begin{defn}
 The constant term of $v$ is obtained by formally setting $t=\log(t)=0$. 
\end{defn}

The justification for this definition is that over the complex numbers
one has $\lim_{t\to 0} t \log(t)=0$. Thus in situations where the
solutions has coefficients in $K[[t]]+tK[[t]][\log(t)] $ this is
indeed the constant term. The same can be argued $p$-adically,
provided one takes the appropriate notion of limit. It is important to
note that the situation described above is indeed what happened for
$p$-adic multiple polylogarithms near $1$, by the main result
of~\cite{F}.

Similarly, there is a basis of (this time global) solutions to $\res_P
\nabla $  with coefficients in $K[\log(\tbar)]  $. In fact,
all the solutions are of the form $\exp(\res_P \Gamma \cdot
\log\tbar)\cdot v $ with $v\in M_P$ and the exponential is a finite sum as
$\res_P \Gamma $ is nilpotent. In this case we can again take the
constant term. This can now be interpreted as simply evaluating at $1$
with the convention that $\log(1)=0$. Thus, we may formally interpret
taking the constant term as continuing to the tangential vector $1$ at
$P$. To make this more than a mere heuristic, though, one needs to
introduce a topology. It can be made precise in the complex
case~\cite{De1} and we will show this also in the $p$-adic case (see
below Proposition~\ref{contan}).

The higher dimensional generalization is now fairly clear.
Suppose that a smooth variety $X$ is given and in
it a divisor $D=\sum_{i\in I} D_i$ with normal crossings. We assume
that all the components $D_i$ are smooth.  For a subset
$J\subset I$ we set $D_J=\cap_{j\in J} D_j$. Let $N_J$ be the normal
bundle to $D_J$ and let $N_J^0$ be the complement in $N_J$ of
$N_{J'}|_{D_J}  $
for $J'\subset J$. Note that one
only need to take $J'$ smaller by one
index and that $N_{\emptyset}$ is considered as the zero
section. Thus,
for example, $N^0_j$ is the (one-dimensional) normal space to $D_j$
minus the zero section. Finally, we denote by $N_J^\zz$ the
restriction of $N^0_J$ to $D_J^0:=D_J- \cup_{j\notin J}D_j$. Following Deligne
we construct, given a
connection on $X$ with logarithmic 
singularities along $D$,  a connection on every
$N_J^\zz$ with logarithmic singularities ``at infinity''. Infinity
here means the union of the hyperplane at infinity for the normal
bundle $N_J$, the hyperplanes $N_J-N_J^0$, and $N_J|_{D_J-D_J^0} $.


For each $j\in J$ consider the valuation $v_j$ on $K(X)$ associated
with the divisor $D_j$. Let $\O_{X}(D^{-1}) $ be the localization of
$\O_X$ at $D$. There exists a  multi-filtration $F_J$ on
$\O_{X}(D^{-1})$, indexed by 
tuples $\chi=(\chi_j)_{j\in J }$, such that $F_J^\chi $ is the
$\O_X$-module generated by
$\{f\in \O_{X}(D^{-1}),\; v_j(f)\ge \chi_j  \text{ for all }  j\in J\}$.
It is easy to see that $\spec (\gr_J \O_{X}(D^{-1})$) is precisely
$N_J^\zz$. Again, a natural map sends a differential form $df$, viewed
as a function on $N_J^\zz$, linear in the bundle direction, to the image
of $f$ in $\gr_J^1  \O_{X}(D^{-1})$.

\begin{rem}
Let $\spec(A)$ be an affine patch on $X$ with coordinates
$t_i\in A$, $i\in I$, such that $D_i\cap \spec(A)$ is defined
by $t_i=0$. Let $B=A/(t_j, j\in J) $, such that $D_J\cap
\spec(A)=\spec(B) $ Let $t^\chi:= \prod_{j\in J} t_j^{\chi_j} $. Then
$F_J^\chi= A t^\chi$ and it follows easily that $\gr_J \cong B[t_j^\pm
, j\not\in J]$
\end{rem}

Now suppose we have a connection with logarithmic singularities along
$D$, $\nabla:M\to M\otimes_{\O_X} \Omega_X^1(\log D)$. We give $
\Omega_{X}^1(D^{-1})=  \Omega_X^1(\log D) \otimes \O_{X}(D^{-1}) $
 the induced filtrations from the filtration on $\O_{X}(D^{-1}) $. It
 is easy to see that the differential $d$ preserves the
 filtration. Now let $M(D^{-1})= M\otimes \O_{X}(D^{-1})$ and
$M\otimes \Omega_{X}^1(D^{-1}) $ have the induced filtrations. It
follows that the extended
connection $\nabla:M(D^{-1})\to M\otimes  \Omega_{X}^1(D^{-1})$
respects the filtration. The connection we
are looking for is the $gr$ of this connection. It is evident from
this construction that if $\nabla $ is flat, so will be the resulting
connection.

\begin{defn}\label{resfunc}
  We call the resulting connection the residue connection on $N_J^\zz$ 
  along $D_J$
  and denote it by $\res_{D,J} \nabla $.
\end{defn}
Note that the construction really depends on both $D$ and $J$ and not
just on $D_J$. It is evident that $\res_{D,J}$ is a tensor functor.

As a consequence of the construction above we obtain more fiber
functors for connections.
\begin{defn}
  If $x\in N_J^\zz $ we let $\omega_{x,J} $ be the fiber functor
  $\omega_{x,J}:= \omega_x \circ \res_{D,J} $.
\end{defn}

Suppose now that $f:(X',D') \to (X,D) $ is a morphism (i.e., takes
$D'$ to $D$) and that $f^\ast D_j = n_j D'_j $. We get an induced map
\begin{equation*}
  N_f:N_{J'}^\zz \to N_{J}^\zz\;, 
\end{equation*}
by looking at the map on the graded sheaf of rings, obtained from $f$,
sending $\gr_J^\chi $ to $\gr_{J'}^{\overline{n}\cdot \chi} $, where
$\overline{n}\cdot \chi $ is the vector with components $n_j\chi_j $.
The following Lemma is trivial.
\begin{lem}\label{makeFrob}
  Suppose that $\Fr$ is the absolute Frobenius. Then $N_{\Fr}$ is also
  the absolute Frobenius. 
\end{lem}

The constructions above can be iterated, and they commute in a sense
we will now explain. For simplicity, we now assume that the divisor
$D$ has only two components, $D_1$ and $D_2$. Consider $N_1$. It has
on it a divisor $D'$ with normal crossings consisting of $D'_1=$ the zero
section and $D'_2=$ restriction of $N_1$ to $D_{12}=D_1\cap D_2$. The
complement of this divisor is exactly $N_1^\zz$. We have the following
result.
\begin{prop}\label{commutant}
  Write $N=N_{D'_1,D'_2}^\zz$ be the construction above performed on
  $N_1$ with respect to $D'$. Then there exists a natural isomorphism
  $\psi: N \xrightarrow{\sim} N_{12}^\zz$ in such a way that for any
  connection $M$ on $X$ with logarithmic singularities along $D$ we
  have a natural isomorphism $\psi^\ast \res_{D,12} M \to
  \res_{D',12}\res_{D,1} M$.
\end{prop}
\begin{proof}
The isomorphism $\psi$ is essentially the isomorphism, existing on any
object with two filtrations, $\gr_{12} \cong \gr_2 \gr_1 $, applied to
the filtrations we have. The isomorphism on the $\res$'s is then a
completely formal consequence of this.
\end{proof}
\begin{rem}\label{Honduras}
To make very concrete the situation of the last Proposition, and to
make it clear how it is going to be used, consider the simplest
situation, $X$ is the affine plane $\bold A^2$ with coordinates $x$ and $y$,
$D=D_1\cup D_2 = \{x=0\} \cup \{y=0\}$. Clearly, $N_1$, the tangent
space to $\{x=0\}$, is again isomorphic to $\bold A^2$, and we may consider
on it the curve $C$ which is the section
$x=1$. Proposition~\ref{commutant}, together with the functoriality of
the $\res$ construction, allows us to get the following: Consider the
fiber functor on connections on $X$ with logarithmic singularities
along $D$: First apply $\res_1$, then restrict to $C$, finally
consider the fiber functor at the tangential point $1$ to $C$ at
$y=0$. Then this fiber functor is naturally isomorphic to the fiber
functor at the normal vector $(1,1)$ at the point $(0,0)$.
\end{rem}

\section{Analytic continuation to tangential points}
\label{sec:analtan}

We now turn to the relation between the tangential theory described in
the previous section and the analytic continuation along Frobenius
explained in
\S \ref{sec:coleman}. Suppose then that $K$ is a $p$-adic field
with ring of integers $\O_K$ and residue field $\kappa$. Let $X$ be a
smooth scheme over
$\O_K$ and let $D=\sum D_i$ be a divisor with relative normal crossings
on $X$. We assume as before that the components $D_i$ are
smooth. The
schemes $N_J^\zz$ are then also smooth over $\O_K$.  Define the
category $\unx$ of unipotent connections on $X_K$
with logarithmic singularities along $D_K$. This is a Tannakian
category over $K$. Suppose that $x\in
(X-D)(\kappa) $ and $y\in
N_J^\zz(\kappa) $ are two points on the special fibers. Restricting
$(M,\nabla)\in \unx$ to $U_x$ and taking horizontal sections we obtain
a fiber
functor $\omega_x$ on $\unx $ as in \S \ref{sec:coleman}. Similarly,
restricting $\res_{D_K,J} \nabla $ to the residue disk $U_y$ of $y$ in
$N_J^\zz$
and taking horizontal sections we obtain a fiber functor
$\omega_y$. Then we have the following result. 
\begin{thm}\label{tanpath}
  For any two points $x$ and $y$, which are either in
  $(X-D)(\kappa) $ or in $ N_J^\zz(\kappa) $ for some $J$,  there exists a
  canonical isomorphism (analytic continuation along a Frobenius
  invariant path),
  $a_{x,y}:\omega_x \xrightarrow{\sim} \omega_y$, such that the
  following properties are satisfied.
  \begin{enumerate}
  \item \label{41a} These isomorphisms are compatible in the sense that
    $a_{y,z} \circ a_{x,y} = a_{x,z}$.
  \item \label{41b} If $x,y\in (X-D)(\kappa)$ then $a_{x,y}$ is the
    isomorphism of
    Theorem~\ref{mainbe} pulled back via $\unx\to \un((X-D)_\kappa)$
  \item \label{41c} If $x,y\in  N_J^\zz(\kappa)$ (for the same $J$)
    then $a_{x,y}$
    is the isomorphism of
    Theorem~\ref{mainbe} pulled back via $\res_{D_K,J}:\unx\to
    \un((N_J^\zz)_\kappa)$
  \end{enumerate}
  The Frobenius invariant path
  is functorial with respect to morphisms $(X',D')\to (X,D) $
  in the obvious sense.
\end{thm}
\begin{proof}
Since the nerve of an affine covering is contractible, we are easily
reduced to the following case: $X=\spec(A)$ is affine and $D$ is
defined on it by the vanishing of $t_1\cdot t_2\cdots t_k$, where
$t_i$ are parameters.
Let $\Adag$ be the weak completion of $A$ in the sense of \cite{Put86} and
let $\Bdag$ be the weak 
completion of $A[(t_1\cdots t_k)^{-1}]$. We clearly have $\Adag
\subset \Bdag$. 
We may further assume that there exists a map, whose reduction is a
Frobenius endomorphism, $\varphi:
\Bdag \to \Bdag $ preserving $\Adag$ and sending $t_i$ to $t_i^q$. We
consider the de Rham complex of $\Adag$ with logarithmic singularities
along $D_K$,
\begin{equation*}
  \dr(\Adag_K,D):\;\Adag_K \xrightarrow{d} \sum \Omega^1(\Adag_K/K)\dlog(t_i)
  \xrightarrow{d} \cdots 
\end{equation*}
and the de Rham complex $\dr(\Bdag_K) $ of $\Bdag_K$. Note that the
cohomology of this last complex is the Monsky-Washnitzer $=$ rigid
cohomology $\hr((X-D)_\kappa)$. There is an obvious
map 
\begin{equation}
  \label{eq:logmap}
 \dr(\Adag_K,D)\to \dr(\Bdag_K)
\end{equation}
and an action of $\varphi^\ast$ on both complexes compatible with this map.
\begin{prop}
  The map  \eqref{eq:logmap} is a quasi-isomorphism.
\end{prop}
\begin{proof}
This follows essentially from the arguments in \cite[Section 6]{Bal-Chi94}
\end{proof}
We only need the fact that the induced  map on $H^1$ is an
injection. This immediately implies the following.
\begin{cor}\label{phiac}
  The map $\varphi^\ast : H^1(\dr(\Adag,D)) \to  H^1(\dr(\Adag,D))$
  is independent of the choice
  of $\varphi$ and has strictly positive weights, 
\end{cor}
\begin{proof}
The same is true for $\hr^1((X-D)_\kappa)$ as in the proof of
Theorem~\ref{mainbe}.
  
\end{proof}

Define now the category $\unad$ of unipotent integrable
$\Adag$-connections with logarithmic singularities along $D$. There is
an action of $\varphi^\ast$ on this category, hence on its fundamental group
and torsors. By the same argument as in the proof of
Theorem~\ref{mainbe} we have using Corollary~\ref{phiac}, that 
there exists a unique invariant path between any two fiber functors on
$\unad$
and that this is independent of the choice of $\varphi$. We have such a
fiber functor $\omega_x$ for $x\in (X-D)_\kappa $, and since  there is an
obvious residual functor $\unad \to \un((N_J^\zz)_\kappa) $, defined
in the same way as in 
Definition~\ref{resfunc}, we also have fiber functors $\omega_y$ for
$y\in(N_J^\zz)_\kappa) $. The fiber functors discussed before the
theorem are the compositions of the ones we just defined with the
restriction functor  $\unx \to \unad$ and we obtain our required path
by pulling back along this restriction. With this definition,
\eqref{41a} and \eqref{41b} are essentially clear, while \eqref{41c}
follows easily because by Lemma~\ref{makeFrob} we have that $\varphi$
induces residually a morphism on $N_J^\zz $  whose reduction is a Frobenius 
automorphism, and consequently the pullback via the
residue functor of a Frobenius invariant path is $\varphi$-invariant.
\end{proof}

We can use the extension of Coleman functions described here to make a
theory of Coleman functions ``of algebraic origin'' on the pair
$(X,D)$. The abstract Coleman functions are then triples $(M,s,h)$
with $M\in \unx$ and $s:M\to \O_{X_K}  $, while $h$ is a compatible system
of local horizontal sections as before. Theorem~\ref{tanpath} allows
us to define a map $\res_{D,J} $ from these algebraic Coleman
functions to Coleman functions on the normal bundle $N_J^\zz$. It is
defined by
\begin{equation*}
\res_{D,J} (M,s,h) = (\res_{D,J} M, \gr s , h)
\end{equation*}
where $h$ now refers to the induced compatible system of horizontal
sections whose existence is guaranteed by Theorem~\ref{tanpath}.

In particular, we obtain from a Coleman function $f$ as above a
Coleman function $f^{(D)} $ on the normal bundle to $D$. The following
propositions suggests how to compute it.

\begin{prop}\label{howcomp}
  Suppose that we have the relation
  $df= \sum \omega_i g_i $ with $f$ and $g_i$ Coleman functions of
  algebraic origin. Then we have the relation $df^{(D)}= \sum (\res_D
  \omega_i) g_i^{(D)}$, where, if
  $\omega$ is locally written as $\omega'+h \dlog(t) $, with $t$ the
  defining parameter for $D$, then $\res_D(\omega)= \omega'|_D + h|_D
  \dlog(\tbar)$.
\end{prop}
\begin{proof}
We consider the way $f$ can be written as a Coleman function following
Remark~\ref{integclos}. In that description suppose that the
connection $M$ constructed there was indeed integrable. Then the
result is easily obtained by computing the residual connection for
$M$. The result remains true if $M$ is not integrable. To see this we
first note that the $M^\textup{int}$ construction is algebraic by
Lemma~2.4 in \cite{Bes00} and the discussion following it. One
further sees that it is a connection with logarithmic singularities
along $D_K$. It follows that we may obtain $f$ from the abstract
Coleman function $(M^\textup{int},s,h) $ as in
Remark~\ref{integclos}. The residual construction is not limited to
integrable connections and is functorial. Thus there is a horizontal map of
connections $\res_{D,i} (M^\textup{int}) \to \res_{D,i} (M)
$. Consequently, $f^{(D)} 
$, which is 
defined on the normal bundle to $D_i$ by $\res_{D,i}
(M^\textup{int},\gr s,h)$, still
satisfies the differential equation given by $\res_{D,i} (M) $
\end{proof}

We can now give, in the $p$-adic case, a tangential basepoint
interpretation of the constant term.
\begin{prop}\label{contan}
  In the situation of Definition~\ref{contstterm} the constant term
  with respect to the parameters $t_1$ and $t_2$ of a horizontal
  section of a unipotent overconvergent isocrystal $\nabla$ on $T$ is the
  evaluation at the normal vector $\tbar_1=1,\tbar_2=1$ of the
  analytic continuation along Frobenius of this horizontal section.
\end{prop}
\begin{proof}
We only prove the $1$-dimensional analogue as the proof is
similar. In fact, we prove something stronger. Suppose (using the
notation of this section) that the
connection is locally given by $d+\Gamma$. Let $f$ be a
horizontal section for the connection $\nabla$ and let $v$ be its
constant term. We claim that the analytic continuation of $f$ is
precisely the solution $c(f):=\exp(\res_P \Gamma \cdot \log\tbar)\cdot v $
(the result is then proved by specializing to $\tbar=1$). We
already know that this is indeed a solution of $\res_P(\nabla)$. The
association $\nabla \mapsto (f\mapsto c(f)) $ is clearly compatible
with direct sums and tensor products hence defines a path. To show
that this is a Frobenius invariant path we may assume that $\varphi$ is
such that $\varphi^\ast(t)=t^p$ The constant term of $\varphi^\ast(f)$
remains $v$ while $\varphi^\ast \nabla$ is associated with $\varphi^\ast \Gamma$. It
follows that 
\begin{equation*}
 c(\varphi^\ast f)=\exp(\res_P (\varphi^\ast \Gamma) \cdot \log\tbar)\cdot v
 =\exp(p \res_P  \Gamma \cdot \log\tbar)\cdot v = c(f)(\tbar^p)\;,
\end{equation*}
which proves our claim.
\end{proof}

\section{The set up of the moduli space $\mathcal{M}_{0,5}$ }\label{review}
This section is devoted to giving a quick review on basic materials of the moduli space
$\mathcal{M}_{0,5}$ of genus $0$ curves with $5$ distinct marked points and
its Deligne-Mumford compactification $\overline{\mathcal{M}_{0,5}}$.
The two (and one) variable $p$-adic multiple polylogarithms are introduced
as Coleman functions on $\mathcal{M}_{0,5}$.

The moduli space $\mathcal{M}_{0,5}=\{(P_i)_{i=1}^5\in(\bold P^1)^5 |
P_i\neq P_j \ (i\neq j)\}/PGL(2)$ is identified with
\begin{equation}\label{coordinate}
\{(x,y)\in\bold A^2\}\backslash 
\{x=0\}\cup\{y=0\}\cup\{x=1\}\cup\{y=1\}\cup\{xy=1\}.
\end{equation}
This identification is given by sending $(x,y)$ to 5 marked points in 
$\bold P^1$ given by $(0,x,1,\frac{1}{y},\infty)$.
The symmetric group $S_5$ acts on $\mathcal{M}_{0,5}$ by
$\sigma(P_i)=P_{\sigma^{-1}(i)}$ ($1\leqslant i\leqslant 5$)
for $\sigma\in S_5$.
Especially for $c=(1,3,5,2,4)\in S_5$ its action is described by 
$x\mapsto\frac{1-y}{1-xy}$, $y\mapsto x$.

The Deligne-Mumford compactification of $\mathcal{M}_{0,5}$ is denoted by
$\overline{\mathcal{M}_{0,5}}$.
This space classifies stable curves of $(0,5)$-type and
the above $S_5$-action extends to the action on $\overline{\mathcal{M}_{0,5}}$.
This space is the blow-up of $(\bold P^1)^2(\supset\mathcal{M}_{0,5})$
at $(x,y)=(1,1)$, $(0,\infty)$ and $(\infty,0)$.
The complement $\overline{\mathcal{M}_{0,5}}-\mathcal{M}_{0,5}$
is a divisor with $10$ components:
$\{x=0\}$, $\{y=0\}$, $\{x=1\}$, $\{y=1\}$, $\{xy=1\}$, $\{x=\infty\}$, 
$\{y=\infty\}$ and $3$ exceptional divisors obtained by blowing up at
$(x,y)=(1,1)$, $(\infty,0)$ and $(0,\infty)$.
In particular for our convenience we denote $\{y=0\}$, $\{x=1\}$, the exceptional divisor at $(1,1)$,
$\{y=1\}$ and $\{x=0\}$ by $D_1$, $D_2$, $D_3$, $D_4$ and 
$D_5$ (or sometimes $D_0$) respectively.
It is because $c^i(D_0)=D_i$.
These five divisors form a pentagon and we denote each vertex
$D_i\cap D_{i-1}$ by $P_i$.
Hence we have $c^i(P_0)=P_i$.
The two dimensional affine space $U_1=Spec \bold Q[x,y]$ gives 
an open affine subset of $\overline{\mathcal{M}_{0,5}}$.
The $S_5$-action gives other open subsets 
$U_i=c^{i-1}(U_1)=Spec \bold Q[z_i, w_i]$ 
($1\leqslant i\leqslant 5$) in $\overline{\mathcal{M}_{0,5}}$
where $(z_1,w_1)=(x,y)$,
$(z_2,w_2)=(y,\frac{1-x}{1-xy})$,
$(z_3,w_3)=(\frac{1-x}{1-xy},1-xy)$,
$(z_4,w_4)=(1-xy,\frac{1-y}{1-xy})$ and
$(z_5,w_5)=(\frac{1-y}{1-xy},x)$.

For ${\bf a}=(a_1,\cdots,a_k)\in\bold Z^k_{>0}$, 
${\bf b}=(b_1,\cdots,b_l)\in\bold Z^l_{>0}$,
and $x$, $y\in\bold Q_p$ with $|x|_p<1$ and $|y|_p<1$
we define {\bf two variable $p$-adic multiple polylogarithm} by
\[
\Li_{\bold a,\bold b}(x,y):=
\underset{<n_1<\cdots<n_l}{\underset{0<m_1<\cdots<m_k}{\sum}}
\frac{\qquad x^{m_k} \qquad y^{n_l}}
{m_1^{a_1}\cdots m_k^{a_k}n_1^{b_1}\cdots n_l^{b_l}}
\in \bold Q_p[[x,y]],
\]
and for ${\bf c}=(c_1,\cdots,c_h)\in\bold Z^h_{>0}$
{\bf one variable $p$-adic multiple polylogarithm} by
\[
\Li_{\bold c}(y):=
{\underset{0<m_1<\cdots<m_h}{\sum}}
\frac{y^{m_h}}{m_1^{c_1}\cdots m_h^{c_h}}
\in \bold Q_p[[y]]\subset
 \bold Q_p[[x,y]].
\]
It is easy to see that these functions satisfy
the following differential equations.
\begin{align}\label{preEngland}
&\frac{d}{dx}\Li_{\bold a,\bold b}(x,y)= \\
&\qquad\quad 
\begin{cases}
\frac{1}{x}\Li_{(a_1,\cdots,a_{k-1},a_k-1),\bold b}(x,y)
\qquad\qquad\qquad\qquad \ \ 
\text{if } a_k\neq 1,\\
\frac{1}{1-x}\Li_{(a_1,\cdots,a_{k-1}),\bold b}(x,y)
-\left(\frac{1}{x}+\frac{1}{1-x}\right)
\Li_{(a_1,\cdots,a_{k-1},b_1),(b_2,\cdots,b_l)}(x,y)\\
\qquad\qquad\qquad\qquad\qquad\qquad\qquad\qquad\qquad\qquad
\text{if } a_k=1,k\neq 1, l\neq 1,\\
\frac{1}{1-x}\Li_{\bold b}(y)
-\left(\frac{1}{x}+\frac{1}{1-x}\right)
\Li_{(b_1),(b_2,\cdots,b_l)}(x,y)
\text{   if } a_k=1,k= 1, l\neq 1,\\
\frac{1}{1-x}\Li_{(a_1,\cdots,a_{k-1}),\bold b}(x,y)
-\left(\frac{1}{x}+\frac{1}{1-x}\right)
\Li_{(a_1,\cdots,a_{k-1},b_1)}(xy)\\
\qquad\qquad\qquad\qquad\qquad\qquad\qquad\qquad\qquad\qquad
\text{if } a_k=1,k\neq 1, l=1,\\
\frac{1}{1-x}\Li_{\bold b}(y)
-\left(\frac{1}{x}+\frac{1}{1-x}\right)
\Li_{\bold b}(xy)
\qquad\qquad\quad \ \ 
\text{if } a_k=1,k=1, l=1,\\
\end{cases}\notag\\
&\frac{d}{dy}\Li_{\bold a,\bold b}(x,y)=
\begin{cases}
\frac{1}{y}\Li_{\bold a,(b_1,\cdots,b_{l-1},b_l-1)}(x,y)
&\text{if } b_l\neq 1,\\
\frac{1}{1-y}\Li_{\bold a,(b_1,\cdots,b_{l-1})}(x,y)
&\text{if } b_l=1, l\neq 1,\\
\frac{1}{1-y}\Li_{\bold a}(xy)
&\text{if } b_l=1,l=1,\\
\end{cases}\notag 
\end{align}

\begin{align}\label{Scottland}
&\frac{d}{dx}\Li_{\bold c}(xy)= 
\begin{cases}
\frac{1}{x}\Li_{(c_1,\cdots,c_{h-1},c_h-1)}(xy)
&\text{if } c_h\neq 1,\\
\frac{y}{1-xy}\Li_{(c_1,\cdots,c_{h-1})}(xy)
&\text{if } c_h=1,h\neq 1, \\
\frac{y}{1-xy}
&\text{if } c_h=1,h=1, \\
\end{cases}\\
&\frac{d}{dy}\Li_{\bold c}(xy)=
\begin{cases}
\frac{1}{y}\Li_{(c_1,\cdots,c_{h-1},c_h-1)}(xy)
&\text{if } c_h\neq 1,\\
\frac{x}{1-xy}\Li_{(c_1,\cdots,c_{h-1})}(xy)
&\text{if } c_h=1,h\neq 1, \\
\frac{x}{1-xy}
&\text{if } c_h=1,h=1, \\
\end{cases}\notag
\end{align}

\begin{align}\label{Ireland}
&\frac{d}{dx}\Li_{\bold c}(y)= 0\\
&\frac{d}{dy}\Li_{\bold c}(y)=
\begin{cases}
\frac{1}{y}\Li_{(c_1,\cdots,c_{h-1},c_h-1)}(y)
&\text{if } c_h\neq 1,\\
\frac{1}{1-y}\Li_{(c_1,\cdots,c_{h-1})}(y)
&\text{if } c_h=1,h\neq 1, \\
\frac{1}{1-y}
&\text{if } c_h=1,h=1, \\
\end{cases}\notag
\end{align}

By the above differential equations,
$Li_{\bold a, \bold b}(x,y)$, $Li_{\bold c}(xy)$ and $Li_{\bold c}(y)$
are all iterated integrals of
$\frac{dx}{x}$, $\frac{dx}{1-x}$, $\frac{dy}{y}$, $\frac{dy}{1-y}$
and $\frac{xdy+ydx}{1-xy}$,
differential forms over $\mathcal{M}_{0,5}$.
Whence they  are obtained from some triple $(M,s,y)$ over
$\mathcal{M}_{0,5}$.
We interpret them as Coleman functions over the rigid triple 
$(\mathcal{M}_{0,5},\overline{\mathcal{M}_{0,5}})$.
This means that they are analytically continued to 
$\mathcal{M}_{0,5}(\bold Q_p)$ as Coleman functions
by the methods of analytically continuation along Frobenius 
in \S\ref{sec:coleman}.

\section{Analytic continuation of two variable $p$-adic multiple polylogarithms}
\label{anacon}

In this section two (and one) variable $p$-adic multiple polylogarithms
are analytically continued into $N_{D_i}^\zz$ ($i\in\bold Z/5$),
which is a Zariski open subset of the normal bundle $N_{D_i}$ 
of the divisor $D_i$ introduced in the previous section.

\begin{nota}
For a Coleman function $f$ over $\mathcal{M}_{0,5}$,
$f^{(D_i)}$ means the analytic continuation of $f$ to $N_{D_i}^\zz$ 
($i\in\bold Z/5$).
For $\bold a=(a_1,\cdots,a_k)\in\bold Z^k_{>0}$ and
$\bold b=(b_1,\cdots,b_l)\in\bold Z^l_{>0}$,
$F_{\bold a,\bold b}$ stands for the Coleman function
$Li_{\bold a,\bold b}(x,y)-Li_{\bold a\bold b}(xy)$ and
for $\bold c=(c_1,\cdots,c_h)\in\bold Z^h_{>0}$,
$G_{\bold c}$ stands for the Coleman function
$Li_{\bold c}(xy)-Li_{\bold c}(y)$ over $\mathcal{M}_{0,5}$.
\end{nota}

\begin{lem}\label{England}
$F^{(D_1)}_{\bold a,\bold b}=0$ and $G^{(D_1)}_{\bold c}=0$ for any index
$\bold a$, $\bold b$ and $\bold c$.
\end{lem}

\begin{pf}
The constant terms of $Li_{\bold a,\bold b}(x,y)$, $Li_{\bold c}(xy)$
and $Li_{\bold c}(y)$ at the origin $P_5$ are zero 
because there are no constant terms in their power series expansions.
We take their differentials and take their residues at $y=0$.
It gives $0$ by induction because each term will be a multiple polylogarithm
with one lower weight than the original one.
Whence $Li_{\bold a,\bold b}(x,y)$, $Li_{\bold c}(xy)$
and $Li_{\bold c}(y)$ are identically zero.
It gives our claim.
\qed
\end{pf}

\begin{lem}\label{France}
$F_{\bold a,\bold b}^{(D_2)}\equiv 0$ if $\bold a$ is admissible
\footnote{
An index $\bold a=(a_1,\cdots,a_k)$ ($a_i\in\bold N$) 
is called {\it admissible} if $a_k>1$.
}
and $G_{\bold c}^{(D_2)}\equiv 0$ for any index $\bold c$.
\end{lem}

\begin{pf}
On the affine coordinate $(z_2,w_2)$ for $U_2$,
the divisor $D_2$ is defined by $w_2=0$.
By using 
$dx=\frac{w_2(1-w_2)}{(z_2w_2-1)^2}dz_2+\frac{z_2-1}{(z_2w_2-1)^2}dw_2$ 
and $dy=dz_2$,
we obtain the following  from
the differential equations in \eqref{preEngland} $\sim$ \eqref{Ireland}.
\begin{align*}
\frac{d}{dz_2}&Li_{\bold a,\bold b}(x,y)=
\begin{cases}
&\frac{w_2}{1-z_2w_2}Li_{(a_1,\cdots,a_{k-1},a_k-1),\bold b}(x,y)+
\frac{1}{z_2}Li_{\bold a,(b_1,\cdots,b_{l-1},b_l-1)}(x,y)\\
&\qquad\qquad\qquad\qquad\qquad\qquad\qquad\qquad
\text{if } a_k>1,b_l\neq 1,\\
&\frac{w_2}{1-z_2w_2}Li_{(a_1,\cdots,a_{k-1},a_k-1),\bold b}(x,y)+
\frac{1}{1-z_2}Li_{\bold a,(b_1,\cdots,b_{l-1})}(x,y)\\
&\qquad\qquad\qquad\qquad\qquad\qquad\qquad\qquad
\text{if } a_k>1,b_l=1, l\neq 1,\\
&\frac{w_2}{1-z_2w_2}Li_{(a_1,\cdots,a_{k-1},a_k-1),\bold b}(x,y)+
\frac{1}{1-z_2}Li_{\bold a}(xy)\\
&\qquad\qquad\qquad\qquad\qquad\qquad\qquad\qquad
\text{if } a_k>1,b_l=1,l=1,\\
\end{cases}\\
\frac{d}{dz_2}&Li_{\bold c}(xy)=
\begin{cases}
\left(\frac{w_2}{1-z_2w_2}+\frac{1}{z_2}\right)\cdot
Li_{(c_1,\cdots,c_{h-1},c_h-1)}(xy)
&\text{if } c_h\neq 1,\\
\left(\frac{z_2w_2(1-w_2)}{(z_2-1)(z_2w_2-1)}+\frac{w_2-1}{z_2-1}\right)\cdot
Li_{(c_1,\cdots,c_{h-1})}(xy)
&\text{if } c_h=1,h\neq 1, \\
\frac{z_2w_2(1-w_2)}{(z_2-1)(z_2w_2-1)}
+\frac{w_2-1}{z_2-1}
&\text{if } c_h=1,h=1, \\
\end{cases}\\
\frac{d}{dz_2}&Li_{\bold c}(y)=
\begin{cases}
\frac{1}{z_2}Li_{(c_1,\cdots,c_{h-1},c_h-1)}(y)
&\text{if } c_h\neq 1,\\
\frac{1}{1-z_2}Li_{(c_1,\cdots,c_{h-1})}(y)
&\text{if } c_h=1,h\neq 1, \\
\frac{1}{1-z_2}
&\text{if } c_h=1,h=1. \\
\end{cases}\\
\end{align*}
Following Proposition~\ref{howcomp}, we compute the residue around $D_2$ 
\begin{align*}
\frac{d}{d\bar z_2}&Li_{\bold a,\bold b}^{(D_2)}(x,y)=
\begin{cases}
\frac{1}{\bar z_2}Li_{\bold a,(b_1,\cdots,b_{l-1},b_l-1)}^{(D_2)}(x,y)
&\text{if } a_k>1,b_l\neq 1,\\
\frac{1}{1-\bar z_2}Li_{\bold a,(b_1,\cdots,b_{l-1})}^{(D_2)}(x,y)
&\text{if } a_k>1,b_l=1, l\neq 1,\\
\frac{1}{1-\bar z_2}Li_{\bold a}^{(D_2)}(xy)
&\text{if } a_k>1,b_l=1,l=1,\\
\end{cases}\\
\frac{d}{d\bar z_2}&Li_{\bold c}^{(D_2)}(xy)=
\begin{cases}
\frac{1}{\bar z_2}Li_{(c_1,\cdots,c_{h-1},c_h-1)}^{(D_2)}(xy)
&\text{if } c_h\neq 1,\\
\frac{1}{1-\bar z_2}
Li_{(c_1,\cdots,c_{h-1})}^{(D_2)}(xy)
&\text{if } c_h=1,h\neq 1, \\
\frac{1}{1-\bar z_2}
&\text{if } c_h=1,h=1, \\
\end{cases}\\
\frac{d}{d\bar z_2}&Li_{\bold c}^{(D_2)}(y)=
\begin{cases}
\frac{1}{\bar z_2}Li_{(c_1,\cdots,c_{h-1},c_h-1)}^{(D_2)}(y)
&\text{if } c_h\neq 1,\\
\frac{1}{1-\bar z_2}Li_{(c_1,\cdots,c_{h-1})}^{(D_2)}(y)
&\text{if } c_h=1,h\neq 1, \\
\frac{1}{1-\bar z_2}
&\text{if } c_h=1,h=1. \\
\end{cases}\\
\end{align*}
Similarly we compute 
\begin{align*}
&\frac{d}{d\bar{w_2}}Li^{(D_2)}_{\bold a,\bold b}(x,y)=0 \qquad\qquad
\text{if  } a_k>1, \\
&\frac{d}{d\bar{w_2}}Li^{(D_2)}_{\bold c}(xy)=0,  \\
&\frac{d}{d\bar{w_2}}Li^{(D_2)}_{\bold c}(y)=0.  \\
\end{align*}
These computations imply that differentials of $F^{(D_2)}_{\bold a \bold b}$
and $G^{(D_2)}_{\bold c}$
with respect to  $\bar{z_2}$ and $\bar{w_2}$ are zero by induction.
Therefore  they must be constant.
By Lemma~\ref{England} their constant terms at $P_1$ is zero.
So they are identically zero.
\qed
\end{pf}

\begin{lem}\label{Russia}
$F_{\bold a,\bold b}^{(D_3)}=0$
if $\bold a$ and $\bold b$ are admissible and
$G_{\bold c}^{(D_3)}=0$
if $\bold c$ is admissible.
\end{lem}

\begin{pf}
On the affine coordinate $(z_3,w_3)$ for $U_3$,
the divisor $D_3$ is defined by $w_3=0$.
By using $dx=-w_3dz_3-z_3dw_3$ and 
$dy=\frac{w_3(1-w_3)}{(z_3w_3-1)^2}dz_3+\frac{z_3-1}{(z_3w_3-1)^2}dw_3$,
we obtain the following  from
the differential equations in \eqref{preEngland} $\sim$ \eqref{Ireland}.
\begin{align*}
\frac{d}{dz_3}&Li_{\bold a,\bold b}(x,y)
=\frac{w_3}{z_3w_3-1}Li_{(a_1,\cdots,a_{k-1},a_k-1),\bold b}(x,y) && \\
\qquad\qquad &+\frac{w_3}{1-z_3w_3}Li_{\bold
  a,(b_1,\cdots,b_{l-1},b_l-1)}(x,y) &&
&\text{ if } a_k>1 \text{ and } b_l>1,\\
\frac{d}{dz_3}&Li_{\bold c}(xy)
=0 &&
&\text{ if } c_h>1,\\
\frac{d}{dz_3}&Li_{\bold c}(y)=
\frac{w_3}{1-z_3w_3}Li_{(c_1,\cdots,c_{h-1},c_h-1)}(y) &&
&\text{ if } c_h>1.
\end{align*}

Following Proposition~\ref{howcomp}, we compute the residue around $D_3$
\begin{align*}
\frac{d}{d\bar z_3}&Li_{\bold a,\bold b}^{(D_3)}(x,y)=0\qquad
\text{ if } a_k>1 \text{ and } b_l>1,\\
\frac{d}{d\bar z_3}&Li_{\bold c}^{(D_3)}(xy)=0\qquad \ \ 
\text{ if } c_h>1,\\
\frac{d}{d\bar z_3}&Li_{\bold c}^{(D_3)}(y)=0\qquad\ \ 
\text{ if } c_h>1.\\
\end{align*}
%
Similarly we also compute
\begin{align*}
\frac{d}{d\bar{w_3}}&Li^{(D_3)}_{\bold a,\bold b}(x,y)=0 \quad \
\text{if  } a_k>1 \text{ and }  b_l>1,\\
\frac{d}{d\bar{w_3}}&Li^{(D_3)}_{\bold c}(xy)=0\qquad
\text{if } c_h>1,  \\
\frac{d}{d\bar{w_3}}&Li^{(D_3)}_{\bold c}(y)=0  \qquad \ \
\text{if } c_h>1.  \\
\end{align*}
These computations imply that differentials of 
$F_{\bold a,\bold b}^{(D_3)}$ and $G_{\bold c}^{(D_3)}$
with respect to $\bar{z_3}$ and $\bar{w_3}$
are zero by induction.
Therefore they must be constant.
By Lemma \ref{France}
their constant term at $P_2$ is zero.
So they are identically zero.
\qed  
\end{pf}

In \cite{F} it was shown that the limit (in a certain way) 
to $z=1$ of $\Li_{k_1,\cdots,k_m}(z)$,
which is a Coleman function over $\bold P^1\backslash \{0,1,\infty\}$,
exists when $k_m>1$ (loc.\ cit.\ Theorem 2.18)
and $p$-adic multiple zeta value $\zeta_p(k_1,\cdots,k_m)$ is defined to be 
this limit value (loc.\ cit.\ Definition 2.17), but
by using the terminologies in \S \ref{sec:tangential}
we reformulate its definition as follows.

\begin{defn}\label{pmzv}
For $k_m>1$, the {\bf $p$-adic multiple zeta value} $\zeta_p(k_1,\cdots,k_m)$
is the constant term of $\Li_{k_1,\cdots,k_m}(z)$ at $z=1$.
\end{defn}

In the case for $k_m=1$, the constant term of 
$\Li_{k_1,\cdots,k_m}(z)$ at $z=1$
is actually equal to the (canonical) regularization
$(-1)^mI_p(BA^{k_{m-1}-1}B\cdots A^{k_1-1}B)$ of
$p$-adic multiple zeta values by loc.\ cit.\ Theorem 2.22
(for this notation, see loc.\ cit.\ Theorem 3.30).\par
The following is important to prove double shuffle relations for
$p$-adic multiple zeta values.

\begin{prop}\label{Germany}
\begin{enumerate}
\item The analytic continuation
$\Li^{(D_3)}_{\bold a,\bold b}(x,y)$
is constant and equal to $\zeta_p(\bold a,\bold b)$
if $\bold a$ and $\bold b$ are admissible.
\item The analytic continuation
$\Li^{(D_3)}_{\bold c}(xy)$ and
$\Li^{(D_3)}_{\bold c}(y)$ are constant and take value 
$\zeta_p(\bold c)$
if $\bold c$ is admissible.
\end{enumerate}
\end{prop}

\begin{pf}
By Lemma \ref{Russia} it is enough to prove this for $Li_{\bold c}^{(D_3)}(y)$.
By the argument in Lemma \ref{England} $Li_{\bold c}^{(D_1)}(y)=0$.
By the computation in Lemma \ref{France} 
$Li_{\bold c}^{(D_2)}(y)=Li_{\bold c}^{(D_2)}(\bar{z_2})$.
So the constant term of $Li_{\bold c}^{(D_2)}(y)$ at $P_2$ 
is equal to the constant term of $Li_{\bold c}(\bar{z_2})$ at $\bar{z_2}=1$,
which is $\zeta_p(\bold c)$.
By the computation in Lemma \ref{Russia}
$Li_{\bold c}^{(D_3)}(y)$ must be constant if $c_h>1$.
Since this constant term must be the constant term of $Li_{\bold c}^{(D_2)}(y)$,
$Li_{\bold c}^{(D_3)}(y)\equiv \zeta_p(\bold c)$ for $c_h>1$.  
\qed
\end{pf}

By discussing on the opposite divisors $D_5$, $D_4$ and $D_3$,
we also obtain the following.

\begin{prop}\label{Italy}
\begin{enumerate}
\item
The analytic continuation
$\Li^{(D_3)}_{\bold a,\bold b}(y,x)$
is constant and equal to $\zeta_p(\bold a,\bold b)$
if $\bold a$ and $\bold b$ are admissible.
\item
The analytic continuation
$\Li^{(D_3)}_{\bold c}(xy)$ and
$\Li^{(D_3)}_{\bold c}(x)$ are constant and equal to
$\zeta_p(\bold c)$
if $\bold c$ is admissible.
\end{enumerate}
\end{prop}

\section{The double shuffle relations}\label{final}
In this section, we prove double shuffle relations for 
$p$-adic multiple zeta values
( Definition \ref{pmzv}).
Firstly we recall double shuffle relations for 
complex multiple zeta values.
Let $\bold a=(a_1,\cdots,a_k)$ and $\bold b=(b_1,\cdots,b_l)$ 
be admissible indices (i.e. $a_k>1$ and $b_l>1$).
The {\bf series shuffle product formulas} 
(called by harmonic product formulas in \cite{F}
and first shuffle relations in \cite{G1})
are relations
\begin{equation}\label{first shuffle}
\zeta(\bold a)\cdot\zeta(\bold b)
={\underset{\sigma\in Sh^{\leqslant}(k,l)}{\sum}}
\zeta(\sigma(\bold a,\bold b))
\end{equation}
which is obtained by expanding the summation on the left hand side
into the summation which give multiple zeta values.
Here 
\begin{align*}
Sh^{\leqslant}(k,l):=\underset{N}{\bigcup}\Bigl\{
\sigma:&\{1,\cdots,k+l\}\to\{1,\cdots,N \}\Bigl| \  \sigma{\text{ is onto}}, \\
&\sigma(1)<\cdots<\sigma(k), \sigma(k+1)<\cdots<\sigma(k+l)
\Bigr\}
\end{align*}
and
$\sigma(\bold a,\bold b)=(c_1,\cdots,c_{N})$ where $N$ is the
cardinality of the image of $\sigma$ and
$$c_i=
\begin{cases}
a_s+b_{t-k} & \text{if  } \sigma^{-1}(i)=\{s,t\} \text{ with } s<t , \\
a_s    & \text{if  } \sigma^{-1}(i)=\{s\}  \quad \text{ with } s\leqslant k.\\
b_{s-k}& \text{if  } \sigma^{-1}(i)=\{s\}  \quad \text{ with } s> k.\\
\end{cases}
$$
One of the easiest example of \eqref{first shuffle} is \eqref{serrel}.

On the other hand, multiple zeta values admit 
an iterated integral expression
(cf. \cite{G1}, \cite{IKZ} see also \cite{Fu})

\begin{align*}
\zeta(\bold a)= \int_0^1 
\underbrace{\frac{du}{u}\circ\dotsm\circ\frac{du}{u}\circ
 \frac{du}{1-u}}_{a_k}& \circ 
\frac{du}{u}\circ\dotsm\dotsm 
\circ\frac{du}{1-u} \\& 
\circ
\underbrace{\frac{du}{u}\circ\dotsm\circ\frac{du}{u}\circ\frac{du}{1-u}}_{a_1}.
\end{align*}
Here for differential $1$-forms 
$\omega_1,\omega_2,\dotsc,\omega_n$ $(n\geqslant 1)$ on $\bold C$
an iterated integral 
$\int_0^1\omega_1\circ\omega_2\circ\dotsi\circ\omega_n$
is defined inductively as
$\int_0^1\omega_1(t_1)\int_0^{t_1}\omega_{2}
\circ\dotsi\circ\omega_n$.
There are the well-known shuffle product formulas 
(for example see ~loc.\ cit.)
of iterated integration 
\[
\int_0^1 \omega_1\circ\cdots\circ\omega_k\cdot
\int_0^1 \omega_{k+1}\circ\cdots\circ\omega_{k+l}=
\sum_{\tau\in Sh(k,l)}
\int_0^1 \omega_{\tau(1)}\circ\cdots\circ\omega_{\tau(k+l)},
\]
where $Sh(k,l)$ is the set of shuffles defined by 
\begin{align*}
Sh(k,l):=\Bigl\{
\tau:&\{1,\cdots,k+l\}\to\{1,\cdots,k+l \}\Bigl| 
\  \tau{\text{ is bijective}}, \\
&\tau(1)<\cdots<\tau(k), \tau(k+1)<\cdots<\tau(k+l)
\Bigr\}.
\end{align*}
They induce the {\bf iterated integral shuffle produce formulas}
(called by shuffle product formulas simply in \cite{F} and
second shuffle relations in \cite{G1})
for multiple zeta values
\begin{equation}\label{second shuffle}
\zeta(\bold a)\cdot\zeta(\bold b)
=\sum_{\tau\in Sh(N_\bold a,N_\bold b)}
\zeta(I_{\tau(W_\bold a,W_\bold b)})
\end{equation}
where $N_{\bold a}=a_1+\cdots+a_k$, $N_{\bold b}=b_1+\cdots+b_l$.
For $\bold c=(c_1,\cdots,c_h)$ with $h,c_1,\dots,c_h\geqslant 1$ the
symbol $W_{\bold c}$ means a word 
$A^{c_h-1}BA^{c_{h-1}-1}B\cdots A^{c_1-1}B$ 
and conversely for given such $W$ we denote its corresponding index by $I_W$.
For words, $W=X_1\cdots X_k$ and $W'=X_{k+1}\cdots X_{k+l}$ with
$X_i\in\{A,B\}$, and $\tau\in Sh(k,l)$ the symbol $\tau(W,W')$
stands for the word $Z_1\cdots Z_{k+l}$ with $Z_i=X_{\tau^{-1}(i)}$. 
One of the easiest example of \eqref{second shuffle} is \eqref{intshuff}.

The {\bf double shuffle relations} for multiple zeta values are
linear relations which are obtained by 
combining two shuffle relations \eqref{double},
i.e.\ series shuffle product formulas \eqref{first shuffle} and 
iterated integral shuffle produce formulas \eqref{second shuffle}
\begin{equation}\label{double}
{\underset{\sigma\in Sh^{\leqslant}(k,l)}{\sum}}
\zeta(\sigma(\bold a,\bold b))
=\sum_{\tau\in Sh(N_\bold a,N_\bold b)}
\zeta(I_{\tau(W_\bold a,W_\bold b)}).
\end{equation}
The following is the easiest example of the double shuffle relations 
obtained from 
\eqref{serrel} and \eqref{intshuff}:
\begin{multline*}
\zeta(k_1,k_2) + \zeta(k_2,k_1) + \zeta(k_1+k_2) \\
=\sum_{i=0}^{k_1-1}\binom{k_2-1+i}{i}\zeta(k_1-i,k_2+i)+\sum_{j=0}^{k_2-1}\binom{k_1-1+j}{j}\zeta(k_2-j,k_1+j)
\end{multline*}
for $k_1,k_2>1$.

\begin{thm}\label{main}
$p$-adic multiple zeta values in convergent case (i.e.\ for admissible indices)
satisfy the series shuffle product formulas, i.e.
\begin{equation}\label{pserrel}
\zeta_p(\bold a)\cdot\zeta_p(\bold b)=
{\underset{\sigma\in Sh^{\leqslant}(k,l)}{\sum}}
\zeta_p(\sigma(\bold a,\bold b))
\end{equation}
for admissible indices $\bold a$ and $\bold b$.
\end{thm}

\begin{pf}
Put $\bold a=(a_1,\cdots,a_k)$ and $\bold b=(b_1,\cdots,b_l)$.
By the power series expansion of 
$\Li_{\bold a,\bold b}(x,y)$ and $\Li_{\bold a}(x)$ in \S \ref{review},
we obtain the following formula
\begin{equation}\label{soleil}
\Li_{\bold a}(x)\cdot \Li_{\bold b}(y)
={\underset{\sigma\in Sh^{\leqslant}(k,l)}{\sum}}
\Li^\sigma_{\bold a,\bold b}(x,y).
\end{equation}
Here 
\[
\Li^\sigma_{\bold a,\bold b}(x,y):=
\underset{(m_1,\cdots,m_k,n_1,\cdots,n_l)\in Z_{++}^\sigma}{\sum}
\frac{\qquad x^{m_k} \qquad y^{n_l}}
{m_1^{a_1}\cdots m_k^{a_k}n_1^{b_1}\cdots n_l^{b_l}}
\]
with
$$
Z_{++}^\sigma=\left\{(c_1,\cdots,c_{k+l})\in\Z_{>0}^{k+l}\bigm|
c_i<c_j \text{ if } \sigma(i)<\sigma(j),
c_i=c_j \text{ if } \sigma(i)=\sigma(j)
\right\}.
$$
Then for each $\sigma\in Sh^{\leqslant}(k,l)$,
$\Li^\sigma_{\bold a,\bold b}(x,y)$ 
can be written 
$\Li_{\bold a^\prime,\bold b^\prime}(x,y)$,
$\Li_{\bold a^\prime,\bold b^\prime}(y,x)$ or
$\Li_{\bold a^\prime,\bold b^\prime}(xy)$
for some indices $\bold a^\prime$ and $\bold b^\prime$.
We note that, if $\bold a$ and $\bold b$ are admissible, then
these $\bold a^\prime$ and $\bold b^\prime$ are also admissible.
By Proposition~\ref{Germany} and Proposition~\ref{Italy},
we know that analytic continuations
$\Li^{(D_3)}_{\bold a,\bold b}(x,y)$,
$\Li^{(D_3)}_{\bold b,\bold a}(y,x)$,
$\Li^{(D_3)}_{\bold a,\bold b}(xy)$,
$\Li^{(D_3)}_{\bold a}(x)$ and
$\Li^{(D_3)}_{\bold b}(y)$
are all constant and take values
$\zeta_p(\bold a,\bold b)$,
$\zeta_p(\bold b,\bold a)$,
$\zeta_p(\bold a,\bold b)$,
$\zeta_p(\bold a)$ and
$\zeta_p(\bold b)$ respectively
when $\bold a$ and $\bold b$ are admissible.
Therefore by taking an analytic continuation along Frobenius 
of both hands sides of \eqref{soleil} into 
$N^{00}_{D_3}(\bold Q_p)$,
we obtain the series shuffle product formulas \eqref{pserrel}
for $p$-adic multiple zeta value in convergent case.
\qed
\end{pf}

By this theorem we say for example
\[
\zeta_p(k_1)\cdot \zeta_p(k_2)= 
\zeta_p(k_1,k_2) + \zeta_p(k_2,k_1) + \zeta_p(k_1+k_2)\; 
\]
for $k_1,k_2>1$
which is a $p$-adic analogue of \eqref{serrel}.

\begin{cor}\label{maincor}
$p$-adic multiple zeta values in convergent case 
satisfy double shuffle relations. Namely
\[
{\underset{\sigma\in Sh^{\leqslant}(k,l)}{\sum}}
\zeta_p(\sigma(\bold a,\bold b))
=\sum_{\tau\in Sh(N_\bold a,N_\bold b)}
\zeta_p(I_{\tau(W_\bold a,W_\bold b)}).
\]
holds for $a_k>1$ and $b_l>1$.
\end{cor}

\begin{pf}
It was shown in \cite{F} Corollary 3.46 that
$p$-adic multiple zeta values satisfy iterated integral 
shuffle product formulas
\begin{equation}\label{gene}
\zeta_p(\bold a)\cdot\zeta_p(\bold b)
=\sum_{\tau\in Sh(N_\bold a,N_\bold b)}
\zeta_p(I_{\tau(W_\bold a,W_\bold b)}).
\end{equation}
By combining it with Theorem~\ref{main}, 
we obtain double shuffle relations for $p$-adic multiple zeta values.
\qed
\end{pf}

Therefore we say for example
\begin{multline*}
\zeta_p(k_1,k_2) + \zeta_p (k_2,k_1) + \zeta_p(k_1+k_2) \\
=\sum_{i=0}^{k_1-1}\binom{k_2-1+i}{i}\zeta_p(k_1-i,k_2+i)+\sum_{j=0}^{k_2-1}\binom{k_1-1+j}{j}\zeta_p(k_2-j,k_1+j)
\end{multline*}
for $k_1,k_2>1$ which is a $p$-adic analogue of \eqref{intshuff}.

In complex case there are two regularizations of multiple zeta values
in divergent case, integral regularization and power series regularization
(see \cite{IKZ}, \cite{G1}\S 2.9 and \S 2.10).
The first ones satisfy iterated integral shuffle product formulas,
the second ones satisfy series shuffle product formulas and
these two regularizations are related by regularization relations.
Actually these provide new type of relations among multiple zeta values.
In the case of $p$-adic multiple zeta values,
$p$-adic analogue of integral regularization 
appear on coefficients of $p$-adic Drinfel'd
associator (see \cite{F}) and they satisfy 
iterated integral shuffle product formulas like \eqref{gene}.
On the other hand, it is not clear at all to say that $p$-adic analogue
of power series regularization satisfy series shuffle product formulas
and regularization relation.
It is because that in the complex case the definition of this regularization 
and the proof of
their series shuffle product formulas and regularization relation
essentially based on the asymptotic
behaviors of power series summations of multiple zeta values
(see \cite{G1} Proposition 2.19)
however in the $p$-adic case our $p$-adic multiple zeta values do not have
power series sum expression like \eqref{defform}.
Recently the validity of these type of relations among 
$p$-adic multiple zeta values were achieved in \cite{FJ}
by using several variable $p$-adic multiple polylogarithm
and a stratification of the moduli $\mathcal{M}_{0,N+3}$
($N\geqslant 3$).

\ifx\undefined\bysame
\newcommand{\bysame}{\leavevmode\hbox to3em{\hrulefill}\,}
\fi

\end{document}